\documentclass{article}

\usepackage{amsmath}

\usepackage{amssymb}

\usepackage{cite}

\usepackage{dsfont}

\usepackage{fancyheadings}

\usepackage{theorem}

\usepackage{verbatim}

\usepackage[all]{xy}

\newcommand{\lieg}{\mathfrak{g}}
\newcommand{\liek}{\mathfrak{k}}
\newcommand{\liea}{\mathfrak{a}}
\newcommand{\lien}{\mathfrak{n}}
\newcommand{\cx}{\mathbb{C}}
\newcommand{\re}{\mathbb{R}}
\newcommand{\sigmalambda}{\sigma_{\lambda}}

\newcommand{\sigmasmlambda}{\sigma^{\ast}_{-\lambda}}
\newtheorem{Def}{Definition}
\newtheorem{Lem}{Lemma}
\newtheorem{Theo}{Theorem}
\newcommand{\VPSL}{V_{\pi^{\sigmalambda}}^{-\omega}}
\newcommand{\VPSLG}{ ^\Gamma V_{\pi^{\sigmalambda}}^{-\omega}}

\newcommand{\VPSsmL}{ V_{\pi^{\sigmasmlambda}}^{-\omega}}
\newcommand{\VPSsmLK}{ V_{\pi^{\sigmasmlambda},K-fin}}
\newcommand{\VPSLK}{V_{\pi^{\sigmalambda}, K-fin}}
\newcommand{\VPL}{V_{\pi^{1_{\lambda}}}^{-\omega}}
\newcommand{\pisl}{\pi^{\sigmalambda}}
\newcommand{\HC}{\mathcal{HC}}
\newcommand{\ZT}{Z}
\title{Analysis on Symmetric and Locally Symmetric Spaces \\
(Multiplicities, Cohomology and Zeta functions)
}
\date{ }
\author{ Ulrich Bunke and Robert Waldm\"uller}
\begin{document}

\maketitle
\section{Multiplicities of principal series}
Let us first fix some notation.
Let $G$ be  a semisimple Liegroup with finite center and
$G\cong KAN$ be the Iwasawa 
decomposition of $G$, where $K$ is a maximal compact subgroup,
$A$ a maximal $\re$-split torus of $G$, and $N$ is a nilpotent subgroup.
Correspondingly, we have a decomposition of the Lie algebra:
$\lieg \cong \liek \oplus \liea \oplus \lien$.
We denote by $P:=MAN$ a minimal parabolic subgroup,
where $M=Z_K(A)$ is the centraliser of $K$ in $A$. For 
$\lambda \in \liea_{\cx}^\ast$, we define
the character $ a \mapsto a^\lambda:=\exp(\lambda(\log a))$.
Furthermore, let $\rho \in \liea^\ast$ denote 
$ H \mapsto \frac{1}{2} Tr (ad(H)|_\lien)$.
For a finite dimensional representation
$(\sigma, V_\sigma)$ of $M$, we define
$\sigmalambda: P \rightarrow GL(V_\sigma)$ by 
$ man \mapsto a^{\rho-\lambda}\sigma(m)$.
\begin{Def} The induced representation
$\pisl:=Ind_P^G(\sigmalambda)$ is called the principal series representation
with parameter $(\sigma,\lambda)$.
\end{Def}
$\pisl$ acts on the space of sections of regularity ? 
of the $G$-equivariant bundle
$V(\sigma,\lambda):=G\times_P V_{\sigmalambda}$. The space of 
sections will be denoted by \\
$V_{\pi_{\sigmalambda}}^{?}:= C^?(G/P,V(\sigma,\lambda))$, where ? belongs to
$K-fin, \omega, \infty,L^2,-\infty$ or $-\omega$.
Here, $K-fin$ means $K$-finite, $\omega$ real-analytic,
$\infty$  smooth and $L^2$  square integrables. 
To be honest,$C^{K-fin}(G/P,V(\sigma,\lambda))$ is only a $(\lieg,K)$-module.
Denote the $\cx$-linear dual of $\sigma$ by $\sigma^\ast$.
To define the hyperfunction and distribution valued sections,
observe that
$V(1,-\rho) \rightarrow G /P $ is the bundle of densities,
and we have a $G$-equivariant pairing
$V(\sigma,\lambda) \otimes V(\sigma^\ast , -\lambda) \rightarrow V(1,-\rho)$.
Combining this pairing with integration, we have a pairing of spaces of 
sections.

\begin{Def}
\begin{eqnarray*}
C^{-\omega}(G/P,V(\sigma,\lambda))& := & 
C^\omega(G/P,V(\sigma^\ast,-\lambda))^\ast \\
C^{-\infty}(G/P,V(\sigma,\lambda))& := & 
C^\infty(G/P,V(\sigma^\ast,-\lambda))^\ast
.
\end{eqnarray*}
\end{Def}
$C^{L^2}(G/P,V(\sigma,\lambda))$ is a $G$-Banach representation.
If $\mathfrak{R}(\lambda)=0$, it is unitary in a natural way. 
\\
\\
\textbf{Example}
We consider $G=SL(2,\re)$.
The Iwasawa components are\\
\begin{center}
$K=\left\{ \left(
\begin{matrix}
 \cos \phi & -\sin \phi \\
\sin \phi & \cos \phi 
\end{matrix}\right) |\phi \in [0,2 \pi) \right\}
\qquad 
A=\left\{ \left(
\begin{matrix}
 t^{\frac{1}{2}} & 0 \\
0 & t^{-\frac{1}{2}} 
\end{matrix}\right) |t \in (0, \infty) \right\}
$
\\
$
N=\left\{ \left(
\begin{matrix}
 1 & x \\
0 & 1 
\end{matrix}\right) |x \in \re \right\}
$
\end{center}
In this case, we have $\liea \cong \re H$, where $H = 
\left(
\begin{matrix}
\frac{1}{2} & 0 \\
0 & -\frac{1}{2}
\end{matrix}
\right)$. Calculating $\rho(H)=\frac{1}{2}$, we obtain that
$\rho \cong \frac{1}{2}.$
Since $M\cong\left\{ \left( 
\begin{matrix}
1 & 0 \\
0 & 1
\end{matrix} \right) , 
\left( \begin{matrix} -1 & 0 \\ 0 & -1 \end{matrix} \right) \right\}
\cong \mathbb{Z}/2\mathbb{Z}  $, $\hat{M}=\{1,\theta\}$ consists of two 
elements. 
Furthermore, $G/P = K/M \cong S^1$. Choosing
$\sigma=1$, the 
principal series representations are parametrized by $\lambda \in \cx$.
For generic $\lambda$, the representation is irreducible. 
\vspace{-4cm}
\\
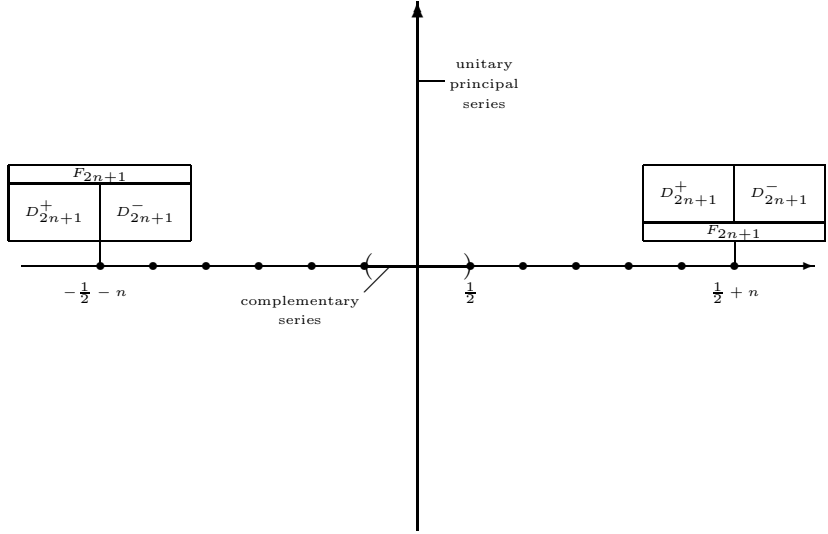
\begin{picture}(200,320)(0,20)
\put(0,100){\vector(1,0){300}}
\put(150,0){\vector(0,1){200}}
\put(170,100){\circle*{3}}
\put(30,100){\circle*{3}}
\put(190,100){\circle*{3}}
\put(50,100){\circle*{3}}
\put(70,100){\circle*{3}}
\put(90,100){\circle*{3}}
\put(130,100){\circle*{3}}
\put(110,100){\circle*{3}}
\put(210,100){\circle*{3}}
\put(250,100){\circle*{3}}
\put(230,100){\circle*{3}}
\put(270,100){\circle*{3}}
\put(28,90){\makebox(0,0){\tiny{$-\frac{1}{2}-n$}}}
\put(170,90){\makebox(0,0){\tiny{$\frac{1}{2}$}}}
\put(270,90){\makebox(0,0){\tiny{$\frac{1}{2} + n$}}}
\thicklines{
\put(130,100){\line(1,0){40}}
\put(150,0){\vector(0,1){200}}
}
\thinlines
\put(133,100){\makebox(0,0){( } }
\put(169,100){\makebox(0,0){)}}
\put(150,170){\line(1,0){10}}
\put(162,160){\tiny{\shortstack{unitary \\ principal \\ series}}}
\put(130,90){\line(1,1){10}}
\put(83,78){\tiny{\shortstack{complementary \\ series}}}
\put(233,138){\tiny{
\begin{tabular}[t]{|c|c|}    \hline
 & \\
$D^+_{2n+1}$ & $D_{2n+1}^- $ \\   
 & \\ \hline
\multicolumn{2}{|c|}{$F_{2n+1}$}\\ \hline
\end{tabular}}}
\put(270,100){\line(0,1){9}}
\put(-7,138){\tiny{
\begin{tabular}[t]{|c|c|}    \hline
\multicolumn{2}{|c|}{$F_{2n+1}$}\\ \hline
 & \\
$D^+_{2n+1}$ & $D_{2n+1}^- $ \\   
 & \\ \hline
\end{tabular}
}}
\put(30,100){\line(0,1){9}}
\end{picture}
\\
\vspace{1cm}
\\
$D^{\stackrel{+}{-}}$ is called the holomorphic/antiholomorphic discrete 
series,
and $F_{2n+1}$ is the $2n+1$-dimensional irreducible representation of 
$SL(2,\re)$. The unitary principal series representations, the 
complementary series representations, $D^{\pm}_{2n+1}$ and $F_1$ are unitary.
Note that $\pi^{1_\pm \frac{1}{2}}$ is a non-trivial extension and therefore
not unitary.
\\
\pagebreak
\\
Let $\Gamma \subset G$ be a cocompact, torsion-free discrete subgroup,
$X:= G / K, \\ Y:= \Gamma \backslash X$ and 
$\VPSLG:=\{\phi \in \VPSL | 
\pi^{\sigmalambda}(\gamma) \phi = \phi \ \ \forall \ \gamma \in \Gamma\},$
the space of $\Gamma$-invariant hyperfunction vectors.
\begin{Lem}[Frobenius reciprocity]
\begin{eqnarray*}
\VPSLG &  \cong &  Hom_G(\VPSsmL,C^{\infty}(\Gamma \backslash G))\\
\phi & \mapsto & 
(f \stackrel{i_{\phi}}{\mapsto} 
(g \mapsto \phi(\pi^{{\sigma^{\ast}}_{-\lambda}}(g)f)))\\
(f \mapsto (\psi(f)(\Gamma 1))) & \leftarrow  & \psi
\end{eqnarray*} 
\end{Lem}
Since $\Gamma$ is cocompact, we have a decomposition
$$C^{\infty}(\Gamma \backslash G) 
\cong \widehat{ \bigoplus}_{(\pi,V_{\pi}) \in \hat{G}_u} 
m_{\Gamma}(\pi)V_{\pi}^{\infty},$$
where $\hat{G}_u$ is the unitary dual of $G$, i.e. the set
of equivalence classes of irreducible unitary representations of $G$.
As a consequence, if $\pi^{\sigmalambda}$ is 
unitary and $\lambda \neq 0$, then it is 
irreducible and therefore $m_{\Gamma}(\pi^{\sigmalambda})=dim(\VPSLG)$.
Moreover, if $\pi^{\sigmalambda}$ is irreducible and not unitary, then
$\VPSLG=\{0\}$.\\
We now return to our example and connect $\VPSLG$ with the eigenvalues of the 
Laplacian on $Y$.
As $K$-homogenous bundles, $V(\sigma, \lambda)=G\times_P V_{\sigmalambda}=
K\times_M V_\sigma$, and for $\sigma = 1$,
$K\times_M V_1 \cong S^{1}\times \cx$. 
Denote by $f \in C^{\infty}(K\times_M V_1)$
the unique $K$-invariant section with $f(1)=1$. For any choice of $\lambda$, 
this corresponds to a real-analytic section $f_{\lambda}$ of 
$V(1,\lambda)$.\\
Using the reciprocity homomorphism 
$i: \VPL \rightarrow 
Hom_G(V^{\omega}_{\pi^{1_{- \lambda}}},C^{\infty}(G))$
for trivial $\Gamma$  , we define the Poisson
transformation
$$i_{\cdot}(f_{- \lambda}):V_{\pi^{1_\lambda}}^{-\omega} \rightarrow 
C^{\infty}(G).$$
Now, any $\phi \in V^{-\omega}_{\pi^{1_{\lambda}}}$ defines
$i_\phi(f_{-\lambda})\in C^{\infty}(X)$.
The Casimir-operator $\Omega$ of $\lieg$ acts on $\pi^{1_{\lambda}}$ 
by a scalar,
$\pi^{1_\lambda}(\Omega)=\mu(\lambda):=\frac{1}{4}-\lambda^2$,
while $\Omega|_{C^\infty(X)}=\Delta_{X}$, the Laplace operator.
Therefore, $i_\phi(f_{-\lambda}) \in Ker (\Delta_{X}-\mu(\lambda))$.
\begin{Theo}[Helgasson \cite{Helg}]\label{TheoH}
For $\lambda \notin \{-\frac{1}{2},-\frac{3}{2},-\frac{5}{2},...\}$,
the Poisson transformation
$$i_{\cdot}(f_{-\lambda}):V_{\pi^{1_\lambda}}^{-\omega} \rightarrow 
Ker(\Delta_{X}-\mu(\lambda))$$
is an isomorphism.
\end{Theo}
To calculate $H^\ast(\Gamma,\VPL)$, we need the 
following classical theorems.
\begin{Theo}\label{TheoA}
Let $M$ be a $C^\omega$-manifold, $E,F \rightarrow M \ 
C^\omega$-vector bundles and \\
$A:C^{\infty}(M,E) \rightarrow C^{\infty}(M,F)$ 
an elliptic operator
with $C^\omega$-coefficients. If $M$ is non-compact, then $A$ is surjective. 
\end{Theo}
\begin{Theo}[\cite{BOgamma},Lemma 2.4, 2.6]\label{TheoB}
Let a discrete group $U$ act properly on a space $M$ and let
$\mathcal{F}$ be a soft or flabby $U$-equivariant sheaf on $M$.
Then $H^i(U,\mathcal{F}(M)) = 0$ for all  $i \geq 1$.
\end{Theo}
Combining Theorem \ref{TheoH} and Theorem \ref{TheoA}, we get
a short exact sequence for 
$\lambda \notin \{-\frac{1}{2},-\frac{3}{2},-\frac{5}{2},...\}$
$$ 0 \rightarrow \VPL \stackrel{i_{\cdot}(f_{-\lambda})}{\rightarrow}
C^{\infty}(X) \stackrel{\Delta_{X}-\mu(\lambda)}{\rightarrow}C^{\infty}(X)
\rightarrow 0.$$
By Theorem \ref{TheoB}, this is a $\Gamma$-acyclic resolution of $\VPL$, so for
$\lambda \notin - \frac{1}{2} - \mathbb{N}_0$, we get
\\
\begin{center}
\begin{tabular}{|c|c|}
\hline
i & $ H^i(\Gamma,\VPL) $\\ \hline
0 & $ker(\Delta_Y - \mu(\lambda))$ \\ \hline
1 & $coker(\Delta_Y - \mu(\lambda))$ \\ \hline
\end{tabular}
\end{center}
Since $(\Delta_Y - \mu(\lambda))$ is an elliptic operator,
both spaces are finite dimensional. To compute
$^\Gamma \VPL$ for $\lambda \in -\frac{1}{2} - \mathbb{N}_0$,
we need $H^i(\Gamma,\VPL)$ for $i \geq 1, \lambda \in \frac{1}{2} + 
\mathbb{N}$ 
as well. \\
In the case of  $\lambda \in \frac{1}{2} + \mathbb{N}$,
we have the exact sequence
$$ 0 \rightarrow F_{2n+1} \rightarrow V^{-\omega}_{\pi^{1_{n+\frac{1}{2}}}} 
\rightarrow
D_{2n+1}^{-\omega} \rightarrow 0$$
Since in this case $\mu(\lambda)=\frac{1}{4}-\lambda^2$ is negative,
$\Delta_Y - \mu(\lambda)$ is a positive operator on a compact manifold and 
therefore injective,so $ ^\Gamma \VPL$ is trivial.
Using the long exact sequence
$$...\rightarrow H^i(\Gamma,F_{2n+1}) \rightarrow 
H^i(\Gamma,V^{-\omega}_{\pi^{1_{n+\frac{1}{2}}}})
\rightarrow  H^i(\Gamma,D_{2n+1}^{\omega}) \rightarrow ...$$
and the topological result \\
\begin{center}
\begin{tabular}{|c|c|c|} \hline
i& $dim(H^i(\Gamma,F_{2n+1}))$ & $dim(H^i(\Gamma,F_1))$ \\ \hline
0 & 0 & 1 \\ \hline
1 & (2g-2)(2n+1) & 2g \\ \hline
2 & 0 & 1 \\ \hline
\end{tabular},\\ 
\end{center}
where $g$ denotes the genus of $Y$, we get
\\
\begin{center}
\begin{tabular}{|c|c|c|c|}\hline
i & $dim(H^i(\Gamma,F_{2n+1}))$ & $dim(H^i(\Gamma,\VPL))$ & 
  $ dim(H^i(\Gamma,D_{2n+1}^{-\omega}))$ \\ \hline
0 & 0 & 0 & (2g-2)(2n+1) \\ \hline
1 & (2g-2)(2n+1) & 0 & 0   \\ \hline
2 & 0 & 0 & 0 \\ \hline
\end{tabular}.
\end{center}
For $\lambda=1/2$, we have that
$dim(ker(\Delta_Y))= dim (coker(\Delta_Y)) = 1$, so we get
\\
\begin{center}
\begin{tabular}{|c|c|c|c|} \hline
i & $dim(H^i(\Gamma,F_1))$ & $dim(H^i(\Gamma,V^{-\omega}_{\pi^{1_{1/2}}}))$ & 
  $ dim(H^i(\Gamma,D_1^{-\omega}))$ \\ \hline
0 & 1 & 1 & 2g \\ \hline
1 & 2g & 1 & 2   \\ \hline
2 & 1 & 0 & 0 \\ \hline
\end{tabular}.
\end{center}
Here, the long exact sequence gave us that 
$dim(H^1(\Gamma,D_1^{-\omega}))$
is either $1$ or $2$ and since $D_1^{-\omega}$ is the sum of two conjugate 
isomorphic
$G$ submodules, it has to be even. 
For $\lambda \in - 1/2 - \mathbb{N}_0$, 
we have a long exact sequence
$$...\rightarrow H^i(\Gamma,D^{-\omega}_{2n+1}) \rightarrow H^i(\Gamma,\VPL)
\rightarrow  H^i(\Gamma,F_{2n+1}) \rightarrow ...$$
Using the result for positive $\lambda$, we obtain immediately for
$\lambda \in -\frac{1}{2} - \mathbb{N}$
\\
\begin{center}
\begin{tabular}{|c|c|c|c|} \hline
i & $dim(H^i(\Gamma,D_{2n+1}^{-\omega}))$ & $dim(H^i(\Gamma,\VPL))$ & 
  $ dim(H^i(\Gamma,F_{2n+1}))$ \\ \hline
0 & (2n+1)(2g-2) &(2n+1)(2g-2)  & 0 \\ \hline
1 & 0 & (2n+1)(2g-2) &  (2n+1)(2g-2)  \\ \hline
2 & 0 & 0 & 0 \\ \hline
\end{tabular}.
\end{center}
In the case of $\lambda = - 1/2$, it is \\
\begin{center}
\begin{tabular}{|c|c|c|c|} \hline
i & $dim(H^i(\Gamma,D_1^{-\omega}))$ & 
$dim(H^i(\Gamma,V_{\pi^{1_{-\frac{1}{2}}}}^{-\omega}))$ & 
  $ dim(H^i(\Gamma,F_1))$ \\ \hline
0 & 2g & 2g & 1 \\ \hline
1 & 2 & 2g+1 & 2g   \\ \hline
2 & 0 & 1 & 1 \\ \hline
\end{tabular}.
\end{center}
Everything except $dim(H^0(\Gamma,V_{\pi^{1_{-\frac{1}{2}}}}^{-\omega}))$ 
is determined by the long exact sequence. For 
$dim(H^0(\Gamma,V_{\pi^{1_{-\frac{1}{2}}}}^{-\omega}))$, we could have 
$2g$ or $2g+1$, but
all invariants are in the submodule $D^{-\omega}_1$ because otherwise we 
would get an embedding of $V_{\pi^{1_{\frac{1}{2}}}}^\omega$ into
$L^2 (\Gamma \backslash G)$ which is impossible, since 
$V_{\pi^{1_{\frac{1}{2}}}}^\omega$ is not unitary.\\ 
In all these cases, the cohomology groups are finite dimensional and
$\mathcal{X}(\Gamma,\VPL)=0$.\\
Let us now consider the general case,
so let $G_{\cx}$ be a connected reductive group over $\cx$,
$G_{\re}$ a real form, $K \subset G_{\re}$ maximal compact,
and $\Gamma \subset G_{\re}$ be cocompact, torsion-free and discrete.
Denote by $Mod(\lieg,K)$ the category of $(\lieg,K)$ modules and by
$\HC(\lieg,K) \subset Mod(\lieg,K)$
the subcategory of Harish-Chandra modules.
Recall that a module is called Harish-Chandra if 
it is finitely generated over the universal enveloping algebra
$\mathcal{U}(\lieg)$ and admissible, i.e. $\forall \gamma \in \hat{K}, \ 
dim(V(\gamma)) < \infty$, where $V(\gamma)$ is the $\gamma$-isotypical 
component of $V$. 
For a $G$ module $V$ denote the $K$-finite submodule by 
$V_{K-fin}=\{v \in V| dim (span(Kv)) < \infty\}$. 
For $V \in \HC(\lieg, K)$ denote by $\tilde{V}:=
Hom(V,\cx)_{K-fin}$ the dual in $\HC(\lieg,K)$.
\begin{Def}
For a Harish-Chandra module $V$, the maximal globalisation of $V$ is
$MG(V):=Hom_{(\lieg,K)}(\tilde{V},C^{\infty}(G)_{K-fin})$.
\end{Def}
Thus, the maximal globalisation $MG(V)$ is a particular continuous 
representation of $G$ whose 
$K$-finite part is isomorphic to $V$.  \\
For a Harish-Chandra module $V$, define the hyperfunction vectors
of $V$ to be $ V^{-\omega}:=(\tilde{V}_B^\omega)^\ast$,
the topological dual of the analytic vectors in any 
Banach globalisation of the dual representation.

\begin{Theo}[Schmid \cite{Sch}]
$MG(V) \cong V^{-\omega}$.
\end{Theo}
\begin{Theo}[Bunke/Olbrich \cite{BOcohom}, 1.4] \label{HGMG}
For $V \in \HC(\lieg,K)$,
$$H^{\ast}(\Gamma,MG(V))\cong \bigoplus_{(\pi,V_{\pi}) \in 
\hat{G_u}}m_{\Gamma}(\pi)
Ext^{\ast}_{(\lieg,K)}(\tilde{V},V_{\pi,K-fin}).$$
\end{Theo}
For the proof of this theorem,
we need the generalisation of Theorem \ref{TheoB}
$$H^i(\Gamma,C^{\infty}(G)_{K-fin}) = 0 \ \forall i \geq 1$$ and
  the following generalisation of Theorem \ref{TheoH}.
\begin{Theo}[Kashiwara-Schmidt \cite{KS}] \label{TheoKS}
$$Ext^{i}_{(\lieg,K)}(\tilde{V},C^{\infty}(G)_{K-fin})=\left\{\begin{array}{c}
MG(V),\ i=0 \\
0,\ i > 0
\end{array} \right.
$$
\end{Theo}
Now, let us chose a projective resolution
$P_\cdot \rightarrow \tilde{V} \rightarrow 0$ of $\tilde{V}$ in $Mod(\lieg,K)$.
By Theorem \ref{TheoKS}, 
$Hom_{(\lieg,K)}(P_{\cdot},C^{\infty}(G)_{K-fin})$ resolves
$MG(V) \ \Gamma$-acyclically ,
therefore $Hom_{\lieg,K}(P_{\cdot},C^{\infty}(\Gamma \backslash G)_{K-fin})$
calculates $H^{\ast}(\Gamma,MG(V))$.
Finally, an analytic argument 
(\cite{BOcohom}, Lemma 3.1)  shows that only a 
finite part of the decomposition
$$C^{\infty}(\Gamma \backslash G)_{K-fin} = 
\widehat{\bigoplus}_{(\pi,V_{\pi}) 
\in \hat{G_u}}
m_{\Gamma}(\pi)V_{\pi,K-fin}^{\infty}$$ contributes to  
$Hom_{(\lieg,K)}(P_{\cdot},C^{\infty}(\Gamma \backslash G)_{K-fin})$.
\section{The Selberg Zeta function}
In this section we will additionally assume that the real rank
of $G$ is equal to one.
A group element $g \in G$ is called hyperbolic if there exist
$m_g \in M, \ a_g \in A$ such that $g$ is conjugate in $G$ to $m_g a_g$ with
$a_g^{\rho} > 1$. If $g$ is hyperbolic, $a_g$ is unique, and $m_g$ is unique 
up to conjugation in $M$.
Note that if $\gamma \in \Gamma, \ \gamma \neq 1$, then $\gamma$ is 
hyperbolic.
\\
For a hyperbolic $g \in G$, we define 
$$ Z(g,\sigma,\lambda):=\prod_{k=0}^{\infty} 
det(1-\sigma_{\lambda-2\rho}(m_g a_g)\otimes S^{k}
(Ad(m_g a_g)|_{\overline{\lien}})),$$
where $\sigma_{\lambda-2\rho}$ denotes the 
representation of $P$ as in section 1 and $S^k$ is the $k'$th symmetric power.
\\
In the example of $G=SL(2,\re)$, a hyperbolic $g$ is conjugate
to $m_g \left( \begin{array}{c c}
t^{1/2} & 0 \\
0 & t^{- 1/2}
\end{array} \right), \\ t > 1$.
Therefore, we get
$$Z(g,1,\lambda) = \prod_{k=0}^{\infty} (1-t^{-\lambda-k - 1/2}).$$
\begin{Def}
The conjugacy class 
$[\gamma] \in C\Gamma$ is called primitive,
if $[\gamma]$ is not of the form
$[\gamma] = [\gamma'^{n}]$ for some $\gamma' \in \Gamma, n > 1$.
\end{Def}
\begin{Def}
The Selberg Zeta function is defined for $\Re(\lambda) > \rho$ by the 
converging product
$$
\ZT(\Gamma,\sigma,\lambda)= \prod_{\begin{array}{c l}
 & [\gamma] \in C \Gamma  \\
 & [\gamma] primitive
\end{array}}
Z(\gamma,\sigma,\lambda).$$
\end{Def}
\begin{Theo}[Fried \cite{Fr}]
$\ZT(\Gamma,\sigma,\lambda)$ has a meromorphic continuation
to $\liea_{\cx}^{\ast}$.
\end{Theo}
The proof uses Ruelles thermodynamic formalism. The disadvantage of this 
method of proof is that is doesn't give any information about the 
singularities of the continuation or about a functional equation.  \\
We use the Selberg trace formuala to calculate the singularities of the 
Zeta function and to determine the functional equation. 
To this end, we introduce  the 
logarithmic derivative
$$L(\Gamma,\sigma,\lambda) := \frac{\ZT'(\Gamma,\sigma,\lambda)}
{\ZT(\Gamma,\sigma,\lambda)}, $$
considered as a 1-form on $\liea_{\cx}^{\ast}$.
Denote by $X^d$ the dual of the symmetric space $X$, i.e. 

$$ 
\begin{array}{|c|c|c|c|c|}
\hline
X & H^m & \cx H^m & \mathbb{H}H^m & Ca H^2 \\ \hline
X^d & S^n & \cx P^m & \mathbb{H} P^m & Ca P^2 \\ \hline
\end{array}.
$$
In the case of $n:=dim X$  odd, we assume that either 
 $\sigma$ is irreducible and the 
isomorphism class of $\sigma$ is 
fixed under the action of the Weyl group, or  
$\sigma = \sigma ' \oplus \sigma '^w$ with $w$ the non-trivial 
element of the Weyl group $W(\lieg,\liea)\cong \mathbb{Z}/2\mathbb{Z}$ and
$\sigma'$ irreducible.\\
The positive root system of $(\lieg,\liea)$ is either of the form 
$\{\alpha\}$ or 
$\{\frac{\alpha}{2},\alpha\}$. We call $\alpha$ the long root and
identify $\liea_{\cx}^{\ast} \cong \cx$ sending $\alpha$ to $1$. 
Define the root vector $H_\alpha \in \liea$ corresponding to a postitive root
$\alpha$ by the condition that
$$\lambda (H_\alpha) = \frac{\langle \lambda,\alpha \rangle}
{\langle \alpha, \alpha \rangle} \quad \forall \lambda \in \liea ^\ast, $$
where $\langle ., .\rangle$ denotes an $Ad(G)$ -invariant inner product 
on $\lieg$. For $\sigma \in\hat{M}$, we define
 $\epsilon_{\alpha}(\sigma) \in \{ 0, \frac{1}{2} \}$ by the 
condition $e^{2\pi i \epsilon_{\alpha}(\sigma)} = 
\sigma(exp(2 \pi i H_\alpha)$ and $\epsilon_\sigma \in \{0, \frac{1}{2}\}$
by requiring $\epsilon_\sigma \equiv |\rho| + \epsilon_{\alpha}(\sigma) \  
\text{mod} \ \mathbb{Z}$.\\ 
\pagebreak

\begin{Theo}[Bunke/Olbrich \cite{BObook}]
There exists a virtual elliptic operator $A_X(\sigma)$ 
and a corresponding operator
$A_{X^d}(\sigma)$, such that in the case of $n$ even 
\begin{eqnarray}
\frac{L(\Gamma,\sigma,\lambda)}{2\lambda}
& = & 'Tr'(\lambda^2 + 
A_{Y}(\sigma))^{-1}  \label{one} \\
&- & \frac{(-1)^{\frac{n}{2}} vol(Y) \pi}
 { vol(X^d)}\frac{P_{\sigma}(\lambda)}{2\lambda} 
\left\{
\begin{array} {l}
tan(\pi \lambda) \  if \ \epsilon_\sigma = \frac{1}{2} \\
-cot(\pi \lambda) \ if \  \epsilon_\sigma = 0
\end{array}
\right.
\label{two}\\
&- & \frac{(-1)^{\frac{n}{2}}vol(Y)}{vol(X^d)} 
 { 'Tr'}  (\lambda^2 - 
A_{X^d}(\sigma))^{-1}\label{three}  
\end{eqnarray}
and in the case of $n$ odd,
\begin{eqnarray*}
\frac{L(\Gamma,\sigma,\lambda)}{2\lambda}
& = & 'Tr'(\lambda^2 + 
A_{Y}(\sigma))^{-1} \\
&+ &   (-1)^{\frac{n+1}{2}} 
\frac{ vol(Y) \pi}
 { vol(X^d)}\frac{P_{\sigma}(\lambda)}{2\lambda} 
\end{eqnarray*}
\end{Theo}
In the theorem,
$P_\sigma$ is some polynomial depending on $\sigma$.
Since $(\lambda^2 + 
A^2_{Y}(\sigma))^{-1}$ is not trace-class,
we have to use a regularised trace $'Tr'$. We refer to
\cite{BObook}, 3.2. for a discussion of the regularisation.\\
>From this formula, it follows that
$L(\Gamma,\sigma,\lambda)$ is meromorphic.\\
In the case of $n$ even, $(\ref{two})$ is odd and $(\ref{three})$ is even. 
The poles of the two terms cancel for $\lambda > 0$, add up for
$\lambda < 0$ and $I:=2 \lambda ((\ref{two}) + (\ref{three}))$ is regular at 
zero 
(\cite{BObook}, 3.2.3). Hence, the set of poles
of $I$ is  $\epsilon_{\sigma}-\mathbb{N}$. \\
For $A$ an operator on some Hilbert space,
$'Tr'(\lambda - A)^{-1}$ has first order poles at the eigenvalues of
$A$ with $res \ (2 \lambda 'Tr'(\lambda^2-A)^{-1})$  equal to the 
dimension of the corresponding eigenspace,
hence 
$$res_{-\epsilon_{\sigma}-k}I =
- \frac{(-1)^{\frac{n}{2}}2 vol(Y)}{vol(X^d)} 
dim E_{A_{X^d}(\sigma)}(-\epsilon_{\sigma}-k)^2.$$
It follows from Hirzebruch proportionality that
$\frac{vol(Y)}{vol(X^d)} \in \mathbb{Z}$ 
(\cite{BObook},Prop. 3.14), so the residues of $I$ are integral. \\
(\ref{one}) has first order poles at $\pm is$ with residue
$dim E_{A_{Y}(\sigma)}(s^2)$, where $s^2$ is a non-zero 
eigenvalue of $A_{Y}(\sigma)$ and if zero is an   
eigenvalue of $A_{Y}(\sigma)$, there is an additional
pole at zero with residue $2 dim ker A_{Y}(\sigma)$. \\
In the case of $n$ odd, the only poles are those coming from the first term.\\
Note that all the residues
of $L$ are integral.\pagebreak\\
\vspace{-3cm}
\\
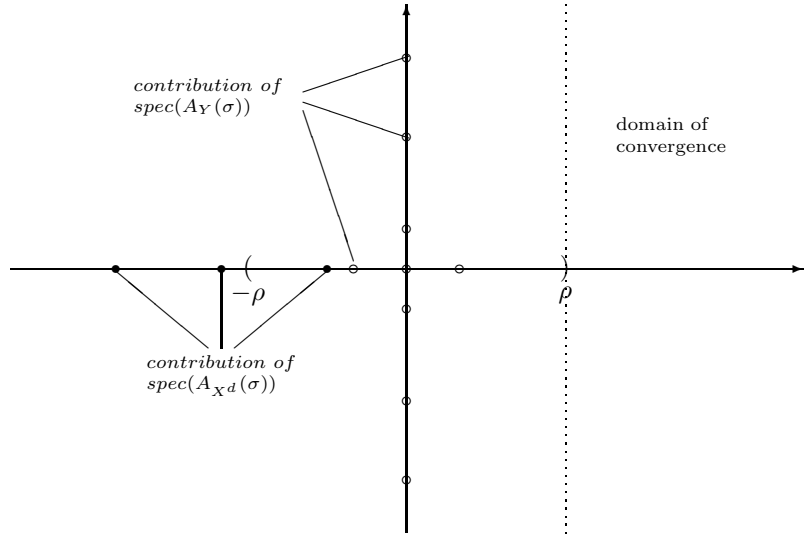
\begin{picture}(200,300)
\put(0,100){\vector(1,0){300}}
\put(150,0){\vector(0,1){200}}
\put(170,100){\circle{3}}
\put(40,100){\circle*{3}}
\put(150,100){\circle{3}}
\put(80,100){\circle*{3}}
\put(130,100){\circle{3}}
\put(120,100){\circle*{3}}
\put(210,100){\makebox(0,0){)}}
\put(210,90){\makebox(0,0){$\rho$}}
\qbezier[65](210,0)(210,100)(210,200)
\put(150,115){\circle{3}}
\put(150,150){\circle{3}}
\put(150,180){\circle{3}}
\put(150,85){\circle{3}}
\put(150,50){\circle{3}}
\put(150,20){\circle{3}}
\put(250,150){\makebox(0,0){\scriptsize{\begin{tabular}{l}
domain of \\ convergence \end{tabular}}}}
\put(111,167){\line(3,1){40}}
\put(111,163){\line(3,-1){40}}
\put(111,160){\line(1,-3){19}}
\put(75,165){\makebox(0,0){\scriptsize{$\begin{array}{l}
contribution \  of \\
spec(A_{Y}(\sigma))
\end{array}$}}}
\put(90,90){\makebox(0,0){$-\rho$}}
\put(90,100){\makebox(0,0){(}}
\put(80,60){\makebox(0,0){\scriptsize{$\begin{array}{l}
contribution \  of \\
spec(A_{X^d}(\sigma))
\end{array}$}}}
\put(80,70){\line(0,1){30}}
\put(75,70){\line(-6,5){35}}
\put(85,70){\line(6,5){35}}
\end{picture}
\vspace{.5cm}
\\
We'll now come back to the example of $G=SL(2,\re), \ \sigma = 1$.
The above theorem specializes to
$$\frac{L(\Gamma,\sigma,\lambda)}{2\lambda} := 
'Tr'(\Delta_{Y}-\frac{1}{4}+\lambda^2)^{-1} + 
\frac{2g -2}{2}{ 'Tr'}(\Delta_{S^2}+\frac{1}{4}-\lambda^2)^{-1} + 
\frac{2g -2}{2}
\frac{\pi \tan(\lambda \pi)}{\lambda},$$
where $g$ is the genus of $\Gamma \backslash X$.
The eigenvalues of $\Delta_{S^2}$ are $\frac{n^2}{4}, \ n \in \mathbb{N}$
with $mult({n^2}/4)=2n+1$. Using this, we get the following picture for 
$ord(\ZT(\Gamma,1,\lambda))$.
\vspace{-3cm}
\\
\begin{picture}(200,300)
\put(120,100){\vector(1,0){190}}
\put(200,0){\vector(0,1){200}}
\put(0,100){\line(1,0){104}}
\put(107,100){\makebox(0,0){\circle*{1}}}
\put(112,100){\makebox(0,0){\circle*{1}}}
\put(117,100){\makebox(0,0){\circle*{1}}}
\put(160,100){\circle*{3}}
\put(20,100){\circle*{3}}
\put(200,100){\circle{2}}
\put(200,100){\circle{3}}
\put(100,100){\circle*{3}}
\put(161,100){\circle{3}}
\put(240,100){\circle{3}}
\put(200,110){\circle{3}}
\put(200,150){\circle{3}}
\put(200,180){\circle{3}}
\put(200,90){\circle{3}}
\put(200,50){\circle{3}}
\put(200,20){\circle{3}}
\put(190,100){\circle{3}}
\put(210,100){\circle{3}}
\put(18,90){\makebox(0,0){\tiny{$-\frac{2n+1}{2}$}}}
\put(18,110){\makebox(0,0){\tiny{$(2g-2)(2n+1)$}}}
\put(158,90){\makebox(0,0){\tiny{$-\frac{1}{2}$}}}
\put(158,110){\makebox(0,0){\tiny{$(2g-2)+1$}}}
\put(240,90){\makebox(0,0){\tiny{$1/2$}}}
\put(240,110){\makebox(0,0){\tiny{$1$}}}
\put(240,178){\makebox(0,0){\tiny{$dim Ker 
 (\Delta_{Y} - \frac{1}{4}+\lambda^2)$}}}
\put(225,175){\line(-1,-1){25}}
\put(225,175){\line(-1,-5){15}}
\put(200,100){\line(1,-2){12}}
\put(249,72){\makebox(0,0)
{\tiny{$2dim Ker(\Delta_{Y} -\frac{1}{4})$}}}
\end{picture}
\\
\begin{Theo}[Bunke/Olbrich, conjectured by Patterson\cite{BOgamma}]
\hspace{7cm}.
For $\lambda \neq 0$ \footnote{For $\lambda =0$ see \cite {BOgamma}},
$$ord_{\lambda} \ZT(\Gamma,\sigma,
\lambda)=-\mathcal{X}'(\Gamma,\VPSL)$$
where $\mathcal{X}'(\Gamma,\VPSL)=\Sigma_{p=0}^{\infty}
(-1)^p p\  dim H^p(\Gamma,\VPSL).$
\end{Theo}
Proof: In the example, the proof is just by inspection. In the general case,
use Theorem \ref{HGMG} to obtain
\begin{eqnarray*}
& & H^p(\Gamma,\VPSL) \\
& = &
\bigoplus_{(\pi,V_{\pi}) \in \hat{G}_u} 
m_\Gamma (\pi) Ext_{(\lieg,K)}^{p}(\VPSsmLK,V_{\pi,K-fin})\\
&= &\bigoplus_{(\pi,V_{\pi}) \in \hat{G}_u}  
m_\Gamma (\pi)Ext_{(\lieg,K)}^{p}(V_{\pi^\ast, K-fin},\VPSLK)\\
&= &\bigoplus_{(\pi,V_{\pi}) \in \hat{G}_u}
m_\Gamma(\pi)Hom_{MA}(H_p(\lien,V_{\pi^\ast})\oplus 
H_{p-1}(\lien,V_{\pi^\ast}),V_{\sigmalambda})\\
&= &\bigoplus_{(\pi,V_{\pi}) \in \hat{G}_u}
m_\Gamma(\pi)Hom_{MA}([H^{n-1-p}(\lien,V_{\pi^\ast, K-fin}) \\
& & \qquad \qquad \qquad \qquad \qquad \qquad 
\oplus H^{n-p}(\lien,V_{\pi^\ast, K-fin})] \otimes 
\Lambda ^{n-1} \lien, V_{\sigmalambda}) \\ 
&= &\bigoplus_{(\pi,V_{\pi}) \in \hat{G}_u}  
m_\Gamma (\pi)\left[((H^p(\lien,V_{\pi,K-fin})\oplus
    H^{p-1}(\lien,V_{\pi,K-fin}))\otimes 
         V_{\sigma^\ast_{-\lambda}})^{MA}\right]^\ast
\end{eqnarray*}
Setting
$$\mathcal{X}(\pi,\sigma,\lambda):=\sum_{p=0}^{\infty} (-1)^p 
dim(H^p(\lien,V_{\pi,K-fin})\otimes V_{\sigma^\ast_{-\lambda}})^{MA},$$
the following theorem finishes the proof.
\begin{Theo}[Juhl, \cite{Juhl}, \cite{Juhlbuch}]
$$ord_\mu \ZT(\Gamma,\sigma,\lambda)=
(-1)^{n-1} \sum_{(\pi,V_{\pi}) \in \hat{G}_u}  
m_{\Gamma}(\pi)\mathcal{X}(\pi,\sigma,\mu).$$
\end{Theo}
\section{Generalisations}
There are two obvious directions in which the theory could be generalised.
First, one could try to generalise the group $G$. Unfortunately,
we do not know of a satisfactory general definition of a Zeta function for 
groups of higher rank, but see \cite{49}, \cite{50}, \cite{Deit} 
for special cases.\\
The other direction is to weaken the conditions on $\Gamma$.
In \cite{BORes} and \cite{BOFuchs}, the theory is generalised
to the case of subgroups $\Gamma \subset G$ of finite covolume. 
Examples of infinite covolume
(convex coccompact) are considered in \cite{BOannals}

\end{document}